\documentclass[letterpaper, 10 pt, conference]{ieeeconf}  

\IEEEoverridecommandlockouts                              

\usepackage[utf8]{inputenc}
\usepackage{bigints}
\usepackage{color,soul}
\usepackage{array}
\newcolumntype{P}[1]{>{\centering\arraybackslash}p{#1}}
\newcolumntype{M}[1]{>{\centering\arraybackslash}m{#1}}
\usepackage[thinlines]{easytable}
\usepackage{multirow}
\usepackage{xcolor,colortbl}

\definecolor{Gray}{gray}{0.85}
\definecolor{LightCyan}{rgb}{0.88,1,1}
\newcolumntype{a}{>{\columncolor{Gray}}c}
\usepackage{hyperref}
\usepackage{amsmath}
\usepackage{caption}
\usepackage{subcaption}
\usepackage{lipsum}
\usepackage{ams}
\usepackage{multicol, blindtext}
\newcommand{\norm}[1]{\left\lVert#1\right\rVert}

   \usepackage[pdftex]{graphicx}     
   \graphicspath{{../pdf/}{../jpeg/}{./image/}}    
   \DeclareGraphicsExtensions{.pdf,.jpeg,.png,.jpg}

\usepackage{amsmath,amssymb,amsfonts}
\usepackage{bm}
\graphicspath{{Figures/}}
\usepackage[normalem]{ulem}

\DeclareMathOperator*{\Minimize}{Minimize:}

\graphicspath{ {./image/} }

\title{\LARGE \bf
Degradation Aware Predictive Energy Management Strategy for Ship Power Systems
}

\author{Satish Vedula, Mehrzad Mohammadi Bijaieh, \textit{Member, IEEE}, Ellis Oti Boateng and \\ Olugbenga Moses Anubi, \textit{Member, IEEE}
\thanks{Authors are with the Department of Electrical and Computer Engineering, Center for Advanced Power systems, Florida State University, Tallahassee, FL 32310, USA
{\tt\small E-mail: \{svedula, mmohammadibijaieh, eotiboateng, oanubi\}@fsu.edu}}%
}

\begin{document}

\maketitle
\thispagestyle{empty}
\pagestyle{empty}

\begin{abstract}
Integration of modern defence weapons into ship power systems poses a challenge in terms of meeting the high ramp rate requirements of those loads. It might be demanding for the generators to meet the ramp rates of these loads. Failure to meet so, might lead to stability issues. This is addressed by conglomeration of generators and energy storage elements to handle the required power demand posed by loads. This paper proposes an energy management strategy based on model predictive control that incorporates the uncertainty in the load prediction. The proposed controller optimally coordinates the power split between the generators and energy storage elements to guarantee that the power demand is met taking into considerations the ramp rate limitations and the load uncertainty. A low bandwidth model consisting of a single generator and a single energy storage element is used to validate the results of the proposed energy management strategy. The results demonstrate the robustness of the controller under load prediction uncertainty and demonstrate the effect of load uncertainty on battery capacity loss.    
\end{abstract}

\section{INTRODUCTION}\label{Sec: Introduction}
Advancements in modern warfare defence mechanisms has led to integration of loads such as electric propulsion motor, radars and other highly non-linear loads into ship power system (SPS). These ever increasing loads could only escalate the ramp rate requirements need to be met by power generation modules (PGMs). In SPS all the electrical machines are considered and configured as a single system namely, integrated power system (IPS). With this configuration it is difficult to accommodate the number of generators required to meet the required load demand adhering to ramp rate requirements. Considering the restrictions on the weight and ship hull dimensions, adding additional generators is not viable. Thus, an alternate energy source with high ramp rate capabilities needs to be added to SPS to mitigate the ramp rate disparities between generators and heavy loads. This problem has been addressed by addition of energy storage systems (ESSs) to the IPS. Thus, the type of distribution system plays a pivotal role in design of IPS components, the transfer from AC to DC distribution system has been proposed in \cite{2007_Doerry}. Most of the existing SPS models consider the DC distribution systems over AC. Since, DC systems are more efficient and easy to design and analyze \cite{2013_Justo}.

Integration of ESSs to SPS provides an additional degree of freedom. ESSs can be used to support the high ramp rate loads. Having multiple ESSs enables seamless operation, as having one ES element out of commission does not hinder SPS operability. The ultimate goal of the SPS operation is to meet the load demand using both generators and ESSs abiding by the ramp rate limitations of respective systems. This task of managing power requirements is steered by an energy management (EM) layer. The main goal of the EM layer is to ensure optimal power distribution subject to prescribed efficiency curves for generators and specific state of charge (SoC) constraints for ESSs, while meeting the load requirements. Most of the existing EM techniques are based on Model predictive control (MPC). 

MPC is an established method used in control of elements involving set of constraints to comply with. MPC was initially designed as a control for slow processes. It was first deployed in chemical process control during the later parts of twentieth century. Substantial developments were made to the MPC control process over the period of time, extending the established concepts to nonlinear MPC and robust MPC. Burgeoning growth in data processing capabilities has led to the use of MPC in fast systems such as power electronics. The power of MPC lies in its ability to operate based on system model information over a prediction horizon and solves the optimization problem at every prescribed time step throughout the horizon. The system limitations such as State of power (SOP) and generator and ESSs ramp rate limitations can be actively handled by the MPC problem. The desired criteria to be achieved by the system is specified as objective function, traditionally as a least squares problem in combination with a cost function associated with generator or ESSs effectiveness. Thus, objective function is also referred by some as cost function \cite{2009_Rawlings}\cite{anubi2015energy}.

\begin{figure*}[h!] 
	\centering
	\includegraphics[width=1.0\textwidth]{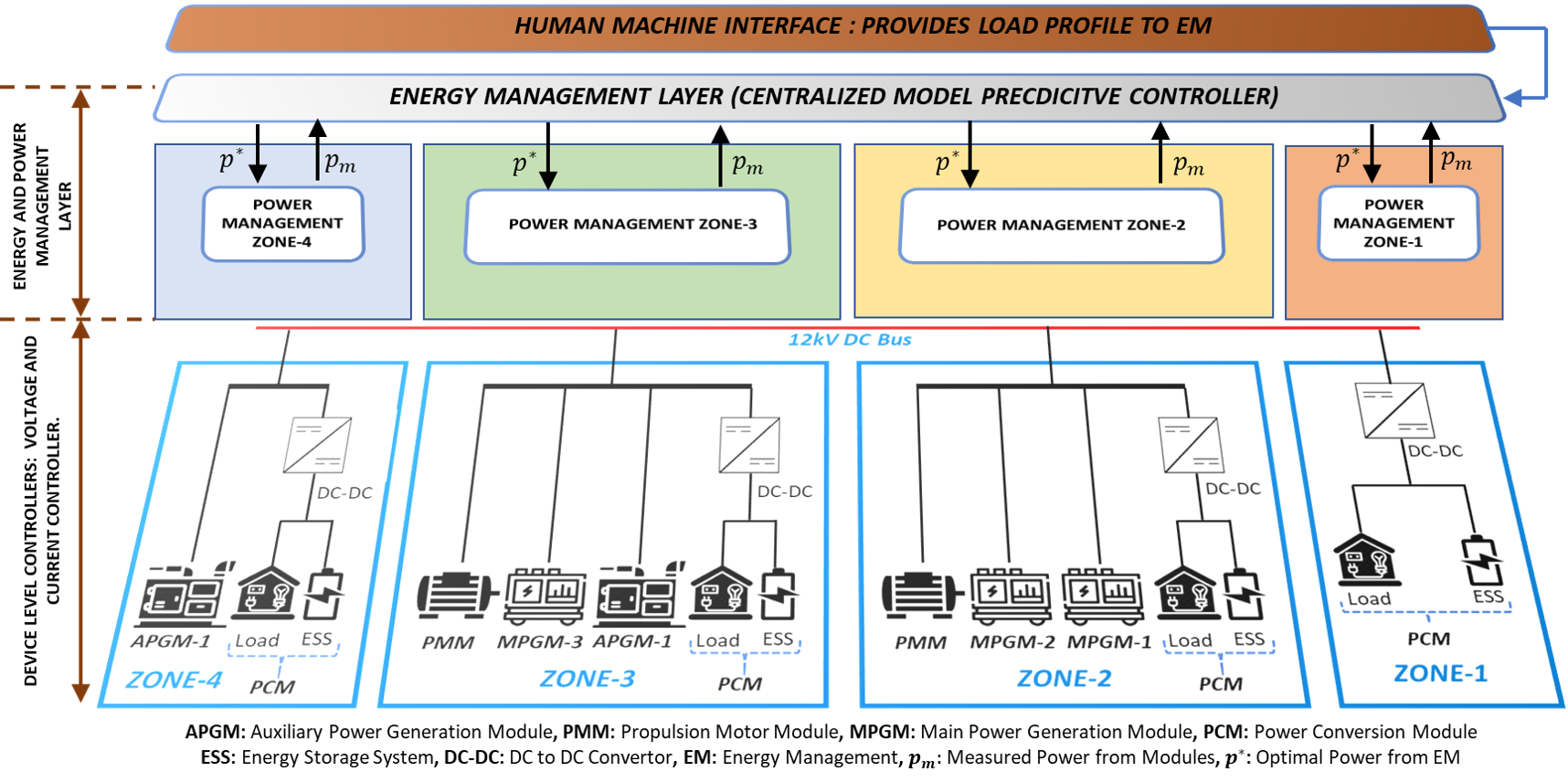} \vspace{-2mm}
	\caption{Energy Management Structure for SPS. 
    }  
	\label{SPS_EM} 
\end{figure*}

MPC based optimization problems can be designed in numerous ways such as: Centralized, Decentralized and Distributed optimization methods. While usage of each of these methods depends on the interest of designer, centralized MPC approaches have been widely used in SPS applications \cite{2017_Gonsoulin}\cite{2017_Vu_1}. While the topic of using centralized or distributed or decentralized MPC for EM and their advantages is widely debated, the concept of distributed control has also been proposed for EM of SPS \cite{2017_Zohrabi}. As discussed before, ramp rate limitations play salient role with high ramp rate loads incorporated into the power system. This calls for an improvised MPC problem which takes into consideration the ramp rate limitations for individual power supply elements in SPS \cite{2017_Vu_3}. Reviewing previously proposed EM methods, this paper presents an MPC energy management method considering the uncertainty in the load prediction. In this work, we have considered the power uncertainty in form of noise injected into load forecast. Varying noise power from 1$\%\ $ to 10$\%\ $ is injected into load forecast and the resulting effect on the robustness of the controller, the battery target SoC and battery capacity loss is studied. The rest of the paper is organized as follows: Section-\ref{Sec: Notation} presents the nomenclature used through rest of the paper, some insight into the SPS power flow model, PCM model used in this work and the battery ESS degradation model considered in this work. In Section-\ref{Sec: Control_Development} the design of the proposed EM strategy is derived. Section-\ref{Sec: Simulation} presents the evaluation of the proposed EM in Simulink environment.

 \section{SYSTEM MODEL}\label{Sec: Notation}
This section presents the mathematical notations used throughout the rest of the paper and introduces basic concepts behind the SPS model such as power flow in SPS, the PCM component model used in this work and the battery capacity loss model.

\subsection{Nomenclature}

This section presents the following notations used throughout the paper: $\mathbb{N}$ denotes the natural number space, $\mathbb{R}$ denotes the real number space, $\mathbb{R}^\text{q}$ denotes real vectors of length $\text{q}$ and $\mathbb{R}^{\text{q}\times\text{p}}$ represents the real matrix with $\text{q}$ rows and $\text{p}$ columns respectively. $\mathbb{R}_+$ represents the positive real number. Lower case alphabets represent the real and natural scalars respectively (e.g $\text{u} \in \mathbb{R}$ and $\text{u} \in \mathbb{N}$) and lower case bold alphabets represent the real vectors (e.g $\textbf{u} \in \mathbb{R}^\text{q}$). The vector of all ones is represented as $\underline{\textbf{1}}$. For a vector $\textbf{u}$,$\textbf{u}_\text{i}$ denotes its $\text{i}^{th}$ element. $\norm{.}$ denotes the two norm.

\subsection{Power-flow Model}
The Ship power system is considered as an islanded MVDC microgrid (MGs). The SPS model presented by the U.S Office of Naval Research (ONR) consists of 4-Zones \cite{ESRDC}. In Fig.\ref{SPS_EM} the bottom half of the picture represents the zonal structure of the SPS. PGMs consists of multiple fuel driven generators and rectifier circuits. The voltage and current controllers which are categorized as device level controllers (DLCs) are also part of PGMs. PCMs consist of power converters, energy storage systems such as batteries and ultra-capacitors and few AC and DC loads. A common 12$\text{kV}$ DC bus connects all the zones. The power generation capacities must meet the power demand requirements under any given circumstance such as: high ramp rate load and load uncertainty. Thus, the net power must always be close to zero or within the acceptable tolerance limits. Thus the SPS power-flow can be represented as: 
\begin{equation}
g(\textbf{p}_\text{g},\textbf{p}_\text{b},\textbf{p}_\text{l}) = 0
\end{equation}
Where $\textbf{p}_\text{g} \in \mathbb{R}^{\text{n}_\text{g}}$ represents the power injected by $\text{n}_\text{g}$ number of generators, $\textbf{p}_\text{b} \in \mathbb{R}^{\text{n}_\text{b}}$ represents the power injected by  $\text{n}_\text{b}$ number of ESSs. $\textbf{p}_\text{l} \in \mathbb{R}^{\text{n}_\text{l}}$ represents $\text{n}_\text{l}$ number of loads in the system. The ESSs can operate bidirectionally, it can be used to support the load requirements and while not in use, it can be operated to reach specified or target SoC.  

\subsection{PCM Model and ESSs Degradation Model}
Lumped ESSs comprising hybrid energy storage system, battery energy storage system (BESS), flywheels is considered as main PCM components in this work. The dynamics of BESS used in this work and the battery current calculations are based on \cite{2020_Bijaieh}. The following dynamics present the relationship between the SoC of the battery and the power injected $\text{p}_\text{b}$:
\begin{equation}\label{SoC_Power}
    \text{SoC} = \frac{\text{Q}_0\text{v}_\text{b}\frac{1}{3600}\int \text{p}_\text{b}(\text{t)\text{dt}}}{\text{Q}_\text{T}\text{v}_\text{b}}
\end{equation}
Where $\text{Q}_0$ is the initial energy stored in BESS in $\text{AHr}$ and $\text{Q}_\text{T}$ is the total energy stored in BESS in $\text{AHr}$. $\text{v}_\text{b}$ represents the bus voltage to which the BESS is coupled.

The battery degradation model considered in this work is based on \cite{2016_Wang}. This model provides the battery capacity loss based on $\text{Ah-throughput}$. The capacity loss is formulated as an exponential function of current throughput. The generator aging or degradation model used is also an exponential model which is a function of power injected \cite{1990_Montanari}.

\section{ENERGY MANAGEMENT DEVELOPMENT} \label{Sec: Control_Development}
The energy management problem is posed as an MPC problem. Consider the following MPC problem:
\begin{subequations}\label{MPC_problem}
\begin{align}
\Minimize \quad & \sum_{\text{k}=1}^{\text{h}}  \norm{\textbf{p}_{\text{g}_\text{k}}+\textbf{p}_{\text{b}_\text{k}}-\textbf{p}_{\text{l}}^\text{f}}^2+\gamma\text{C}(\textbf{p}_{\text{k}})\\
\textrm{s.t.} \quad &  \sum_{\text{k}=1}^{\text{h}}\textbf{p}_{\text{b}_\text{k}}=  \dfrac{3600\text{Q}_\text{T}\text{v}_\text{b}^*}{\text{T}_\text{s}}(\text{q}_0-\text{q}_\text{h})\\
& \textbf{p}_{\text{g}_0} = \textbf{p}_{\text{g}_\text{m}} \\
& \textbf{p}_{\text{b}_0} = \textbf{p}_{\text{b}_\text{m}} \\
& \underline{\textbf{p}}_\text{g} \preceq \textbf{p}_{\text{g}_\text{k}} \preceq \overline{\textbf{p}}_\text{g}  \\
& \underline{\textbf{p}}_\text{b} \preceq \textbf{p}_{\text{b}_\text{k}} \preceq \overline{\textbf{p}}_\text{b} \\
& | \textbf{p}_{\text{g}_\text{k}}-\textbf{p}_{\text{g}_\text{k-1}}| \preceq \text{r}_\text{g}, \hspace{3mm} \text{k} = 1,2,....,\text{h} \\
& | \textbf{p}_{\text{b}_\text{k}}-\textbf{p}_{\text{b}_\text{k-1}}| \preceq \text{r}_\text{b}, \hspace{3mm} \text{k} = 1,2,....,\text{h}
\end{align}
\end{subequations}

The initialization powers for the MPC problem are $\textbf{p}_{\text{g}_0}$ and $\textbf{p}_{\text{b}_0}$. These represent the instantaneous power measurements coming into the optimizer from the system generators and ESSs. $\textbf{p}_{\text{g}_\text{m}}$ and $\textbf{p}_{\text{b}_\text{m}} $ represent the instantaneous measurements of generator power and ESSs power. Where $\textbf{p}_{\text{g}_\text{k}}$ is the power injected by the generators at the given instant $\text{k}$. $\textbf{p}_{\text{b}_\text{k}}$ is the power injected or absorbed by the ESSs. $\textbf{p}_\text{l}^\text{f}$ is the load forecast with uncertainty.The noise is added as uncertainty to the forecasted load. $\text{h} \in \mathbb{N}$ is the horizon of the MPC problem. The problem in (\ref{MPC_problem}a) is solved subject to the constraints over the specified horizon $\text{h}$ at every time instant $\text{k}$. $\text{C}: \mathbb{R}\longrightarrow\mathbb{R}_+$ is the cost capturing the attributes such as efficiency for generators. For ESSs the cost can be associated with effective battery degradation. $\gamma \geq 0$ is a positive scalar weight which is used to address the effect of cost in the objective function. The case of $\gamma = 0$ represents that the cost is not considered in the objective problem. The equation (\ref{MPC_problem}b) represents the equality constraint of the MPC problem. $\text{Q}_\text{T}$ represents the total capacity of the ES element. $\text{v}_\text{b}^*$ represents the measured or instantaneous voltage at the bus to which the ESS is coupled. $\text{q}_0$ is the SoC measured at the beginning of the horizon, this comes as an input from the system as a measurement to the controller. $\text{q}_\text{h}$ is the desired target SoC at the end of the horizon. $\text{T}_\text{s}$ is the time step of the simulation. Equations (\ref{MPC_problem}e) and (\ref{MPC_problem}f) represent the box constraints of the optimization problem. $\underline{\textbf{p}}_\text{b}$ represents the lower power limitations on the ESS and $\overline{\textbf{p}}_\text{b}$ represents the upper power restrictions on the ESS. $\underline{\textbf{p}}_\text{g}$ represent the lower power limitations on the generator power and $\overline{\textbf{p}}_\text{g}$ represent the upper power limitations on the generator power. While $\underline{\textbf{p}}_\text{b}$ can be negative due to bidirectional nature of the ESS, $\underline{\textbf{p}}_\text{g}$ is strictly greater than zero or equal to zero. Equations (\ref{MPC_problem}g) and (\ref{MPC_problem}h) represent the inequality constraints associated with the MPC problem. $\text{r}_\text{g}$ corresponds to the ramp rate limitations of the generators and $\text{r}_\text{b}$ corresponds to the ramp rate limitations of the ESSs. The MPC problem mentioned in (\ref{MPC_problem}a)-(\ref{MPC_problem}h) is hard to use for simulation purpose. Thus, a simplified version of it is derived which is easy to implement using existing solvers such as fmincon and quadprog in MATLAB-Simulink. The reformulation of the MPC problem is as follows:
\begin{subequations}\label{MPC_reformulation}
\begin{align}
\Minimize \quad & \frac{1}{2}\text{x}^\text{T}\text{H}\text{x}+\text{f}^\text{T}\text{x}\\
\textrm{s.t.} \quad &  \text{A}_\text{eq}\textbf{x} = \text{b}_\text{eq}\\
& \text{A}\textbf{x} \preceq \textbf{b}  \\
& \textbf{LB} \preceq \textbf{x} \preceq \textbf{UB} 
\end{align}
\end{subequations}

Where $\mathbf{x} = \begin{bmatrix} \textbf{p}_{\text{g}_\text{k}} \\ \textbf{p}_{\text{b}_\text{k}} \end{bmatrix}$ represents the vector of generator and battery powers. The above reformulation is a quadratic program. The conversion to the above reformulation from the initially proposed MPC form in (\ref{MPC_problem}) is presented in \cite{2011_Jerez}.

\section{SIMULATION EXAMPLE}\label{Sec: Simulation}
\begin{figure}[h!] 
	\centering
	\includegraphics[width=0.48\textwidth]{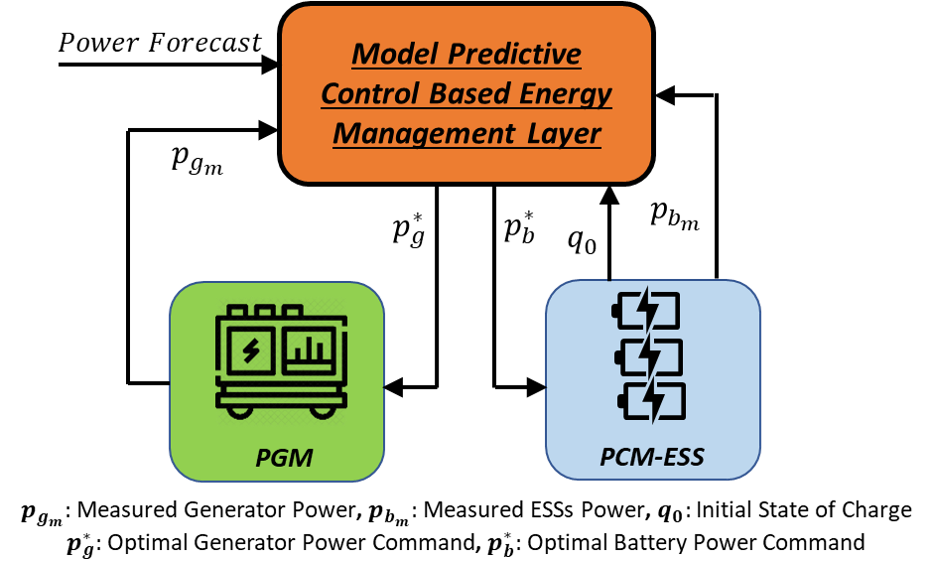} 
	\caption{Simplified Single PGM, PCM and Load Model Used for Simulation. 
    }
	\label{SPS_Sim}
\end{figure}

The simplified single PGM, PCM and Load model shown in Fig.\ref{SPS_Sim} is considered for simulation purpose. The PGM, PCM and the Load models used in this work are low fidelity mathematical models based on \cite{bijaieh2021model}. The simplified model is implemented in the Simulink environment. Table-\ref{sys_rating} shows the component ratings such as generational capacities, ramp rate limitations, lower and upper power limitations for both generator and battery elements that have been used in this simulation. The simulation is based on the problem in (\ref{MPC_problem}a)-(\ref{MPC_problem}h). The overall EM optimizer and the PGM and PCM models including the DLCs are designed and implemented in Simulink. The time step considered for the MPC problem is $10^{-4}$sec. The rate transition of $10^{-3}$sec is considered between the EM layer and the system, accommodating for the time constants of the device level components of the system. The effect cost of generation in the objective function in this work is considered to be zero i.e $\gamma$ is considered to be zero. 

\begin{table}[h!]
\renewcommand{\arraystretch}{1.3}
\centering
\caption{System Components Power Ratings} \vspace{-1mm}
\label{sys_rating}
\begin{tabular}{lcccccccc}
\hline
\hline
\rowcolor{Gray}
\multicolumn{9}{c}{\textbf{PGM, PCM and Load Ratings}}                                                 \\ \hline 
               & \multicolumn{2}{c}{Power (MW)}  & \multicolumn{2}{c}{RR (MW/sec)}  & \multicolumn{2}{c}{LL (MW)} & \multicolumn{2}{c}{UL (MW)}  \\ \hline
PGM    & \multicolumn{2}{c}{29}  & \multicolumn{2}{c}{2.9}  & \multicolumn{2}{c}{0.29} & \multicolumn{2}{c}{27.5}  \\ 
  \hline
PCM    & \multicolumn{2}{c}{30}  & \multicolumn{2}{c}{10}  & \multicolumn{2}{c}{-10.64} & \multicolumn{2}{c}{10.64}  \\ \hline
 
Load    & \multicolumn{2}{c}{30}  & \multicolumn{2}{c}{10}  & \multicolumn{2}{c}{-} & \multicolumn{2}{c}{-}  \\ \hline \hline

\end{tabular}
\end{table}

\begin{figure}[h!] 
	\centering
	\includegraphics[width=0.48\textwidth]{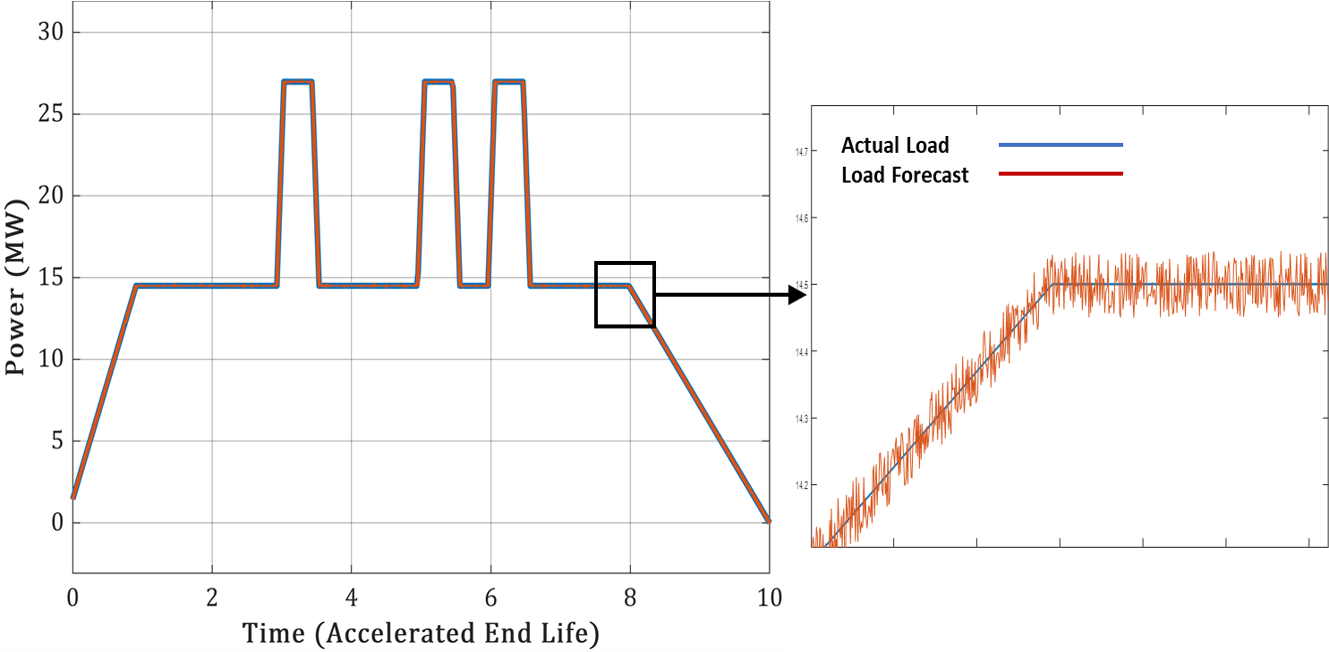} 
	\caption{Actual Load vs the Lower forecast with uncertainty. 
    }
	\label{Power_Forecast}
\end{figure}

\begin{figure}[h!] 
	\centering
	\includegraphics[width=0.48\textwidth]{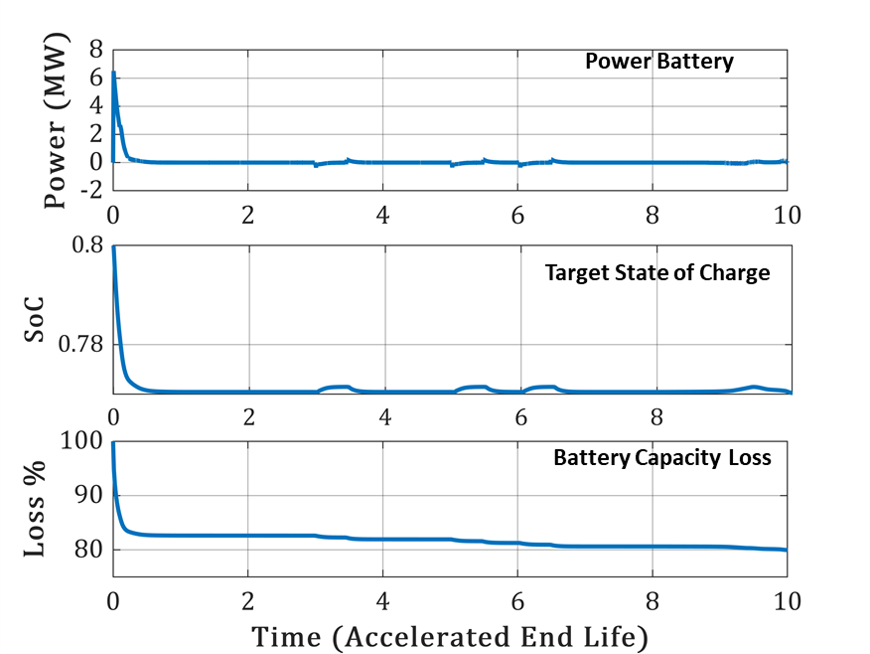} 
	\caption{Optimization power input to the battery, the state of charge of the battery in response to the power injection and the battery capacity loss. 
    }
	\label{SoC_Sim}
\end{figure}

\begin{figure}[h!] 
	\centering
	\includegraphics[width=0.45\textwidth]{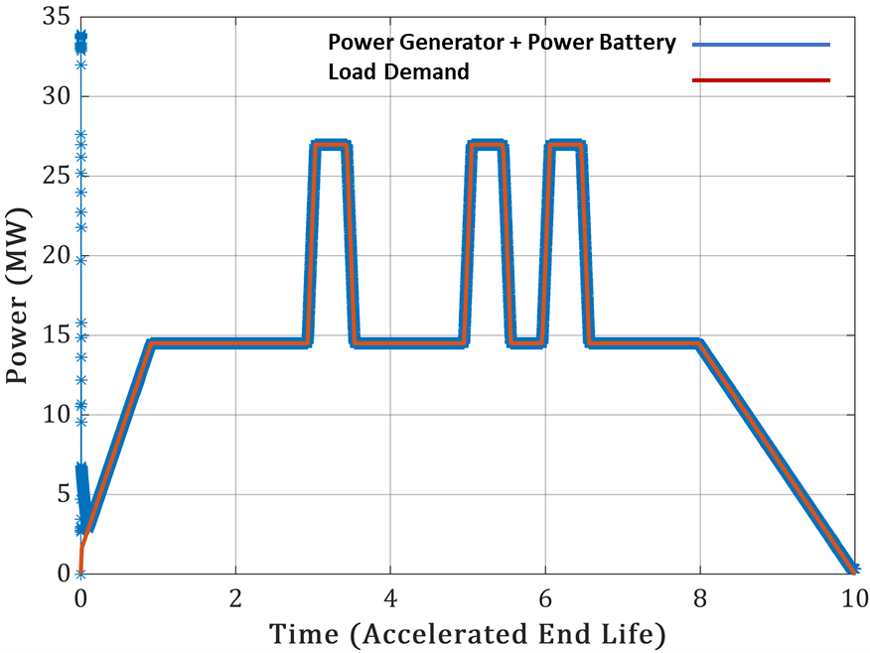} 
	\caption{Closed Loop Simulation result for power injected by generator and battery to meet the required load demand. 
    }
	\label{Closed_Loop}
\end{figure}

The following scenario has been implemented and their results have been presented: the scenario demonstrates the behavior of controller under load forecast with $10\%\ $ added forecast uncertainty. The control power input to the battery from the optimizer and the variations in the SoC of the battery and the battery capacity loss $\%\ $ are also studied and presented for the scenario mentioned above. In the above mentioned scenario the target SoC $\text{q}_\text{h}$ is assumed to be constant. The target SoC or the SoC at the end of the horizon is considered as 0.77. The initial SoC which acts as an input to the optimizer from the system can be seen from the Fig.\ref{SoC_Sim} as 0.8. The SoC of the battery and the battery capacity loss are also shown in Fig.\ref{SoC_Sim}. Fig.\ref{Power_Forecast} shows the actual load profile superimposed on the load with 10$\%\ $ uncertainty. The power forecast is generated with assumption of deploying pulsed power loads (PPL) multiple times in a given duration. The peaks in Fig.\ref{Power_Forecast} represent PPLs. The simulation performed is a closed loop MPC problem. The simulation run-time is 10 seconds. The degradation of components is accelerated in order to accommodate for the time constraint and also reach the component end of life. Thus, the time axis for figures refers to time as \textit{accelerated end life}. Due to the initialization values of generator and ESSs, initially the simulation starts at an infeasible point as seen in Fig.\ref{Closed_Loop}. But, as the controller starts enforcing the load forecast on the system, the system decision variables generator power and battery power begin converging to a feasible point and start tracking the desired load power.

\begin{figure}[h!] 
	\centering
	\includegraphics[width=0.45\textwidth]{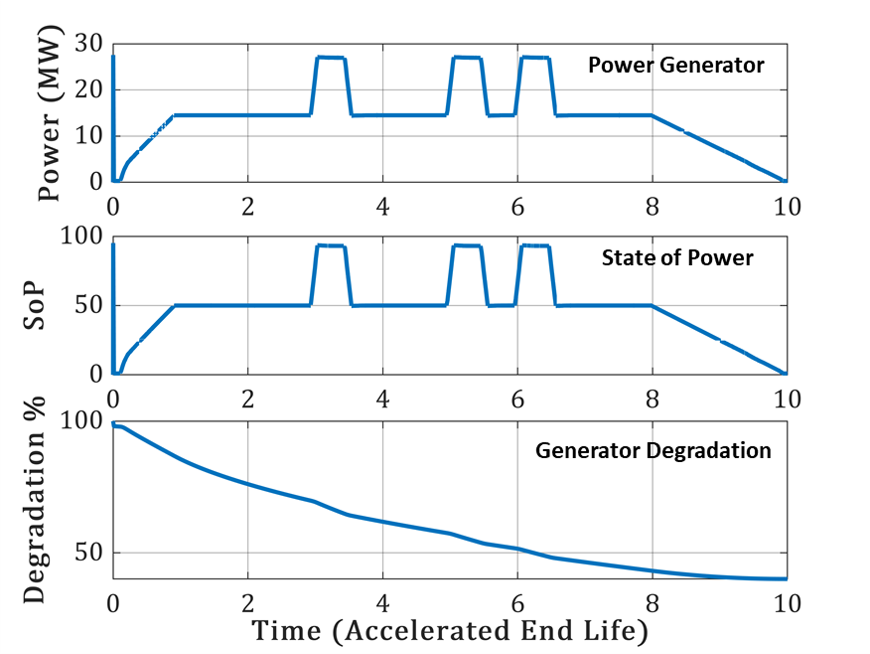} 
	\caption{Optimization power input to the generator, the state of power of the generator in response to the power injection and the generator degradation. 
    }
	\label{SoP_Sim}
\end{figure}

\begin{figure}[h!] 
	\centering
	\includegraphics[width=0.48\textwidth]{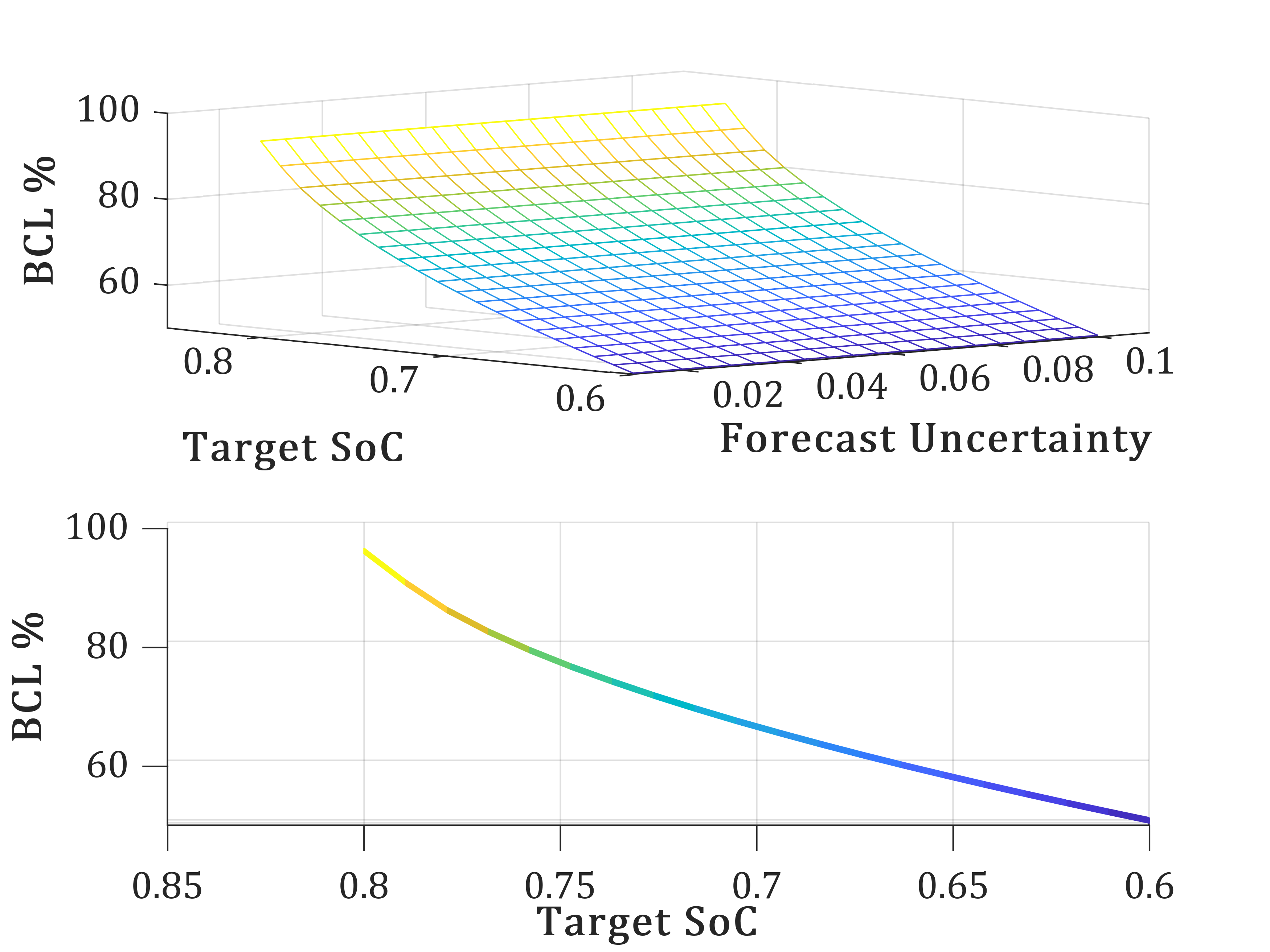} 
	\caption{\textbf{BCL}: Battery Capacity Loss. First figure shows the relationship between capacity loss, target SoCs and forecast uncertainty. Second figure shows the quadratic relation between capacity loss and target SoCs. 
    }
	\label{Noise_Sim}
\end{figure}

The generator power, the generator state of power (SoP) and the generator degradation are plotted in Fig.\ref{SoP_Sim}. In this work the generator SoP is considered to be the instantaneous power of generator divided by the rated power of the generator. The interrelation between noise power intensity, target SoC and the battery capacity loss is presented in Fig.\ref{Noise_Sim}. The noise power used in this case is the noise sweep performed from power intensities 1$\%\ $ to 10$\%\ $ and in this case the target SoC is swept from 0.8 to 0.6. Parallel simulation \textit{parsim} environment in simulink is used to simulate the system multiple times in parallel with above mentioned noise power and target SoC sweeps. The data on capacity loss for a certain target SOC and forecast uncertainty is obtained from 400 parallel simulations and is plotted in Fig.\ref{Noise_Sim}. The simulation results shown in Fig.\ref{Noise_Sim}, the \textit{degradation curve}, demonstrate the robustness of controller under load uncertainty. The simulations show the relationship between target SoC and the battery capacity loss to be quadratic. 

 


\section{CONCLUSIONS}\label{Sec: Conclusion}

This paper presents a model predictive control based energy management strategy for SPS with load forecast uncertainty. The uncertainty is introduced as noise power and its effect on controller robustness, battery capacity loss are studied. Different noise power simulations are implemented in parallel in simulink. The simulations demonstrate the robustness of the controller. The simulations also show the effect of forecast uncertainty on battery capacity and the rate of BESS usage. The target SoCs role in battery degradation is also presented. Finally, the degradation curve for ESS is presented. While, the proposed the controller is implemented on a simple PGM, PCM and Load model, the future objective is to extend this work to the 2-zone and 4-zone SPS model and add a battery degradation aware term into the objective function as a decision variable which solves the optimal power split problem based on current ESS degradation. The idea is to make use of the obtained degradation curve as battery cost. Thus, the future iteration of this work will include a predictive control problem which takes into consideration the battery degradation as cost function.

\section{ACKNOWLEDGEMENT}\label{Sec: ACKNOWLEDGEMENT}
This material is based upon research supported by, or in part by, the U.S. Office of Naval Research under award number N00014-16-1-2956. Any opinions, findings and conclusions or recommendations expressed in this material are those of the authors and do not necessarily reflect the views of the ONR.







\bibliography{FinalVersion/References_Updated.bib}
\bibliographystyle{ieeetr}

\end{document}